\documentclass{amsart}
\usepackage[]{amsmath}
\newcommand{\myhref}[2]{#2}

\newcommand{\ext}[2]{#1/#2}
\newcommand{\Gal}[2]{\mathop\mathrm{Gal}(\ext{#1}{#2})}
\newcommand{\kro}[2]{\left( \frac{#1}{#2} \right)}
\newcommand{\Z}{{\bf Z}}
\newcommand{\Q}{{\bf Q}}
\newcommand{\CEPL}{\mathop{\mbox{}L}\nolimits}

{
\theoremstyle{definition}
\newtheorem{definition}{Definition}
}

\newtheorem{theorem}{Theorem}

\begin{document}

\title{First-hit analysis of algorithms for computing quadratic irregularity}
\author{Joshua Holden}
\address{Department of Mathematics,
Duke University,
Durham, NC 27708, USA}
\email{holden@math.duke.edu}
\urladdr{http://www.math.duke.edu/\~{}holden}

\keywords{Bernoulli numbers, irregular primes, zeta functions,
quadratic extensions, cyclotomic extensions, class groups,
computational number theory, cryptography}


\subjclass{Primary 11Y40, 11Y60, 11Y16, 11R42;
  Secondary 11B68, 11R29, 94A60, 11R18}

\begin{abstract}
    The author has previously extended the theory of regular and
    irregular primes to the setting of arbitrary totally real number
    fields.
It has been conjectured that the Bernoulli numbers, or
    alternatively the values of the Riemann zeta function at odd
    negative integers, are evenly distributed modulo $p$ for every
    $p$.  This is the basis of a well-known heuristic, given by Siegel
    in~\cite{Siegel64},
    for estimating the frequency of irregular primes.
    So far, analyses have shown
    that if $\Q(\sqrt{D})$ is a real quadratic field, then the values
    of the zeta function $\zeta_{D}(1-2m)=\zeta_{\Q(\sqrt{D})}(1-2m)$
    at negative odd integers are also distributed as expected modulo $p$
    for any $p$.
We use this heuristic to predict the computational time required to
find quadratic analogues of irregular primes with a given order of
magnitude.  We also discuss alternative ways of collecting large
amounts of data to test the heuristic.
\end{abstract}

\maketitle

\section{Introduction}
\label{introduction}
Let  $\Q(\sqrt{D})$ be a real quadratic field with $D$ a positive fundamental
discriminant.  In several previous papers the author has defined an
analogue for the theory of regular and irregular primes in this
setting, based on the following definition:

\begin{definition}
Let $\zeta_{D}$ be the zeta function for $\Q(\sqrt{D})$, and let
$\delta$ be equal to $p-1$ unless $D = p$, in which case $\delta =
(p-1)/2$.  We say that $p$ is \emph{$D$-regular} if $p$ is relatively
prime to $\zeta_{D}(1-2m)$ for all integers $m$ such that $2 \leq 2m
\leq \delta - 2$ and also $p$ is relatively prime to $p \zeta_{D}(1 -
\delta)$.  The number of such zeta-values that are divisible by $p$
will be the \emph{index of $D$-irregularity} of $p$.
\end{definition}

(More generally, we
may refer to the concept as ``quadratic irregularity'';
see~\cite{Holden98, Holden99, mathcomp} for more details and
extensions to any totally real number field.)

According to a well-known
theorem of Kummer, $p$ divides the order of the class group of
$\Q(\zeta_{p})$ if and only if $p$ divides the numerator of a
Bernoulli number $B_{2m}$ for some even $2m$ such that $2 \leq 2m \leq
p-3$.  Such primes are called irregular; the others are called
regular.
In~\cite{Holden99}, building on work of Greenberg and Kudo, the author proved that
in the setting we have described above Kummer's criterion can be extended
to give information about whether $p$ divides the class number (that
is, the order of the class group) of
$\Q(\sqrt{D}, \zeta_{p})$. To be exact, we have:

\begin{theorem}[Greenberg, Holden] \label{GH}
Assume that $p$ does not divide $D$.
Then $p$ divides the class number
of $\Q(\sqrt{D}, \zeta_{p})$ if and only if $p$
is not $D$-regular.
\end{theorem}

The main focus of this paper is in finding large $p$ which are
irregular for some $D$.  This may be useful for cryptography, in that
 one common way of constructing
public-key cryptographic systems is to utilize the problem of finding
a discrete logarithm in some abelian group.
In order to make sure that the discrete logarithm
problem is computationally hard, one needs to know something about the
structure of the group involved, e.g. that it is divisible by a large
prime.  Theorem~\ref{GH} shows that if $p$ is a large $D$-irregular prime
and $p$ does not divide $D$, then the class group of $\Q(\sqrt{D}, \zeta_{p})$
may be suitable for cryptography.  We will come back to this in
Section~\ref{practical}.

\section{Search Algorithms}
\label{algorithms}
Suppose, for instance, that we want to find $p$ of a specified size
dividing \mbox{$\zeta_{D}(1-2m)$} for some $m$ such that $2 \leq 2m
\leq \delta-2$ or dividing $p\zeta_{D}(1-\delta)$.
(For reasons that
will become clear, we will not encounter the situation $D=p$ in
practice, so we may focus on the case where $\delta = p-1$.  See
Section~\ref{probabilities} for the details.)
More
specifically, we might fix a real number $c$ greater than $1$ and then look for
$m$ and $D$ such that $P \leq p \leq cP$ and $p$ divides
$\zeta_{D}(1-2m)$ for some positive $m$ less than or equal to $(p-1)/2$.  (In practice, $c=2$
would probably be the most common choice, since that would be
equivalent to specifying the size of $p$ in number of bits.)

The algorithm that we will use to carry out this search is described
in~\cite{mathcomp}.  The algorithms there fall into two basic types.
The first type calculates $\zeta_{D}(1-2m)$ in a range of $m$ for each
$D$ before going on to the next $D$, and calculates each value in time
$O(m^{O(1)} D^{1+o(1)})$ when amortized over both $D$ and $m$.  The
second type calculates $\zeta_{D}(1-2m)$ in a range of $D$ for each
$m$ before going on to the next $m$.  If one keeps a table of
intermediate values as described in Section~3 of~\cite{mathcomp}, this
algorithm can calculate each value in time $O(m^{O(1)}
\CEPL(D)^{O(1)})$ when amortized over both $D$ and $m$, where
$\CEPL(x)$ is a subexponential function corresponding to a choice of
factoring routine used in the calculation, e.g.
$\CEPL(x)=e^{c(\log x)^{1/3} (\log \log x)^{2/3}}$ for the number field
sieve.  The facts that the amortized times are subpolynomial in $D$
and that a range of $D$ are calculated for each $m$ suggest using this
algorithm.

In fact, one might suppose that one could always use $m=1$ and search
until an appropriate $D$ is found without ever going on to the next
$m$.  However, one other factor needs to be taken into account.  The
size of the numerator of \mbox{$\zeta_{D}(1-2m)$} is $O(m (\lg m + \lg
D))$ bits (see~\cite{Holden98}), so \mbox{$\zeta_{D}(1-2m)$} is much
more likely to have large prime factors for large $m$ than for small
$m$.  The same issue comes up for $D$, of course, but to a lesser degree
and in a way which does not greatly affect this algorithm, since a
large number of values of $D$ are used for each $m$.

To make sure that we can avoid getting stuck in a range where the
values of the numerator of \mbox{$\zeta_{D}(1-2m)$} are too small, we
will give our algorithm parameters $M_{1}$,  $D_{1}$, and
$D_{2}$ such that we always have $m \geq M_{1}$ and  $D_{1} \leq D \leq D_{2}$.
In Section~\ref{hypotheses} we will explain some conjectures
which imply that for each pair $(D, m)$,
the chance that each prime between $P$ and
$cP$ divides \mbox{$\zeta_{D}(1-2m)$} is approximately $2/P$.  Given
this, the Prime Number Theorem implies that the chance that some prime
between $P$ and $cP$  divides \mbox{$\zeta_{D}(1-2m)$} is approximately  $2(c-1) /
(\log P)$.

Now if we are trying $D_{1} \leq D \leq D_{2}$ for each $m$, we see
that the chance that we find a suitable $D$ for any given $m$ is approximately
$$\left(\frac{3}{\pi^{2}}\right)\frac{2(D_{2}-D_{1})(c-1)}{\log P},$$
since the asymptotic density of fundamental discriminants in the
integers is $3/\pi^{2}$.
However, if we consider the typical case $c=2$, then we see that we
only have to choose $D_{2}-D_{1}$ to be on the order of magnitude of $\log
P$ for the expected probability of success on any given $m$ to be $1$!
Thus with $D_{2}-D_{1}$ reasonably large, the expected time using
this strategy is
$$O\left(\left( \frac{\log P}{c-1} \right) \CEPL(D)^{O(1)} + P\right),$$
where the added term of $P$ accounts for the time it takes to check
whether each $p$ divides each computed \mbox{$\zeta_{D}(1-2m)$}.
Table~\ref{first-hit-timing} provides some actual timing examples of
this algorithm, running on a Pentium III computer using the  Linux
operating system and the GP-Pari interpreted language.
(See~\cite{PARImanual}.) In all cases $D_{1}=5$ and $M_{1}=2$.

\begin{table}
\caption{Time in minutes to find the first suitable pair $(D,m)$ with given
parameters}
    \label{first-hit-timing}

    \begin{tabular}{|c|c|c|c|c|c|c|}
\hline
$c$ & $P$ &  $D_{2}$ & $(D,m)$ & $p$ & min.\\
\hline
1.01 & $10^{5}$  & $10^{4}$ & $(4156,2)$ & 100391 & 27 \\
1.1 & $10^{5}$  & $10^{4}$ & $(697,2)$ & 106681 & 33 \\
1.6 & $10^{5}$ & $10^{4}$ & $(205,2)$ & 113173 & 81 \\
2 & $10^{5}$ & $10^{4}$ & $(184,2)$ & 164999 & 82 \\
\hline
2 & $10^{6}$ &  $300$ & $(40,3)$ & 1264807 & 169 \\
2 & $10^{6}$ &  $10^{3}$ & $(380,2)$ & 1017299 & 191 \\
2 & $10^{6}$ &  $10^{4}$ & $(380,2)$ & 1017299 & 191 \\
\hline
2 & $2\cdot 10^{6}$ & $10^{4}$ & $(317, 2)$ & 2027569 & 569 \\
\hline

\end{tabular}
\end{table}

\section{The Hypotheses: Conjectures and Previous Results}
\label{hypotheses}
\label{olddata}
The hypotheses mentioned in Section~\ref{algorithms} stem from the
conjecture, made (not very explicitly) by Siegel in~\cite{Siegel64},
that the numerators of the
Bernoulli numbers $B_{2m}$ were evenly distributed modulo $p$ for any
odd prime $p$.  Siegel used the conjecture to derive a conjectural
density for the irregular primes.  (Lehmer seems to have done the
same thing in~\cite{Lehmer} but only gives the density.)
Siegel's hypothesis was used more generally by Johnson (\cite{Johnson75})
and independently by Wooldridge (\cite[Chap.~III]{Wooldridge}) to predict the
density of primes with a given index of irregularity, that is such
that $p$ divides a given number of the Bernoulli numbers $B_{2},
\ldots, B_{p-3}$.
It also comes in handy for predicting many other
values that are related to irregular primes, such as the order of
magnitude of the first prime of a given index of irregularity.  (See,
for example, \cite{Wagstaff78}.)
Since $B_{2m}=-\zeta(1-2m)(2m)$, it is equivalent to say that the
values of $\zeta(1-2m)$ are evenly distributed modulo $p$, where
$\zeta(s)$ is the Riemann zeta function.

Little or no progress has been made on proving Siegel's hypothesis,
but a great deal of data has been collected, especially in regard to
the prediction of Johnson and Wooldridge.  Specifically, this
prediction says that as $p\to\infty$, the probability that $p$ has
index of irregularity $r$ goes to
$$\left( \frac{1}{2} \right)^{r} \frac{e^{-1/2}}{r!} \enspace .$$
(In addition to the original sources, the details may be found in
Section~5.3 of~\cite{Washington}.)  Note that this prediction does
not rely on the full strength of Siegel's hypothesis, but merely on
the weaker hypothesis that the Bernoulli numbers
are $0$ modulo $p$ with probability $1/p$.  The assumptions made in
Section~\ref{algorithms}
relate only to predictions about indices of irregularity
based on this weaker hypothesis.

Wagstaff, in
\cite{Wagstaff78}, computed $u_{r}(x)$, the fraction of primes not
exceeding $x$ with index $r$ of irregularity for each $r$ between 0
and 2 and for all $r\geq 3$ grouped together, and compared this
distribution to the predicted distribution for each multiple $x$ of
1000 up to 125000.  The result of the chi-squared test ``fluctuated
usually between 0.1 and 1.0 and had the value 0.29 at $x=125000$.  It
was 0.03 at $x=8000$'' \cite{Wagstaff78}.  These results correspond to
significance levels of .992, .801, .962, and .999, respectively.
(The significance levels used in this paper correspond roughly to the
probability that the agreement between the observed results and the
predicted results is \emph{not} due to chance.  Statisticians
consider the threshold for considering a result to be not due to
chance to be a significance level of .9 to .95.  Since we are not
actually conducting a valid statistical study in this paper, all of
the statistical results should be taken with a very large grain of
salt.)

Buhler, Crandall, Ernvall, and Mets\"ankyl\"a hold the record for
computations with irregular primes, having found all the irregular
primes below four million as described in \cite{BCEM}.  They do not
seem to have done a chi-squared analysis, but they tabulate the
values of $u_{r}(x)$ for $x=4000000$ and $r$ between 0 and 7.
A chi-squared test using the same methodology as before has the
result 1.02, for a significance level of .796.
Earlier, in \cite{BCS}, Buhler, Crandall, and Sompolski tabulated the
same data for $x=1000000$.  The result of the same chi-squared test is
0.78, for a significance level of .854.

Unfortunately, the
only way to collect data to test Siegel's hypothesis is to
investigate $B_{2m}$ for larger and larger $m$, which is very
computationally intensive.
(See~\cite{Bach}
or~\cite{Fillebrown} for details.)

However, in the more general number field case, there are many more
dimensions to the problem.  We start by restricting our attention to
the case of $k$ an abelian totally real number field.  Then we know that
$$\zeta_{k}(s)=\prod_{\chi\in\hat{G}} L(s, \chi)$$
where $\hat{G}$ is the character group of $G=\Gal{k}{\Q}$ and $L(s,
\chi)$ is the $L$-function associated with the character $\chi$.
Note that $L(s, 1)=\zeta(s)$, so the Riemann zeta function is a factor
of the zeta function for $k$.
(See~\cite{CF}, e.g., for more details.)
Certainly it seems likely that for a fixed (totally real) number field
$k$ and character $\chi$ the values of the numerator of $L(1-2m,
\chi)$ are evenly distributed modulo $p$ as $m$ varies.  (It is known
that these values are rational numbers.)
We also hypothesize that these values for different $\chi$ are
independent, which implies that the numerators of $\zeta_{k}(1-2m)$
are distributed modulo $p$ like the product of $\left| G \right|$
independent integer variables, each of which is evenly distributed
modulo $p$.  We will refer to this as the ``product distribution'',
for lack of a better term.
However, it also is reasonable to conjecture
that for a fixed $m$ the values of $\zeta_{k}(1-2m)$ are distributed
according to the product distribution modulo $p$ as
$k$ varies.  More precisely, if we fix $m$
and the degree of $k$ we expect the values to be distributed
according to the product distribution modulo $p$ as
the discriminant of $k$ varies.  Alternatively, if we fix $m$ and the
discriminant of $k$ we expect the values to be distributed according
to the product distribution modulo $p$ as
the degree varies.

In this paper we will be considering the former situation.  As in the
previous sections, we fix the
degree at $2$, and let $k=\Q(\sqrt{D})$ be a real quadratic field
with zeta function $\zeta_{D}(s)$.
In this case
$$\zeta_{D}(s)= L(s, 1) L(s, \chi) = \zeta(s) L(s, \chi)$$
where $\chi(s) = \kro{D}{s}$, the Kronecker symbol, where appropriate.

In addition to the above definitions we will make one more set:

\begin{definition}
Let $\chi$ be as above and let $\delta$ be as in
Section~\ref{introduction}.  We will say that $p$ is \emph{$\chi$-regular}
if $p$ is relatively prime to $L(1-2m, \chi)$ for all integers $m$
such that $2 \leq 2m \leq \delta - 2$ and also $p$ is relatively prime
to $p L(1-2m, \chi)$.  The number of such $L$-values that are
divisible by $p$ will be the \emph{index of $\chi$-irregularity} of
$p$.
\end{definition}

Saying that the values of $\zeta_{k}(1-2m)$ are distributed according
to the product distribution and that the values of $\zeta(1-2m)$ are
evenly distributed is the same as saying that the values of $L(1-2m,
\chi)$ are evenly distributed modulo $p$.  Then we can make the same
prediction about the indices of $\chi$-irregularity that Johnson and
Wooldridge made about the indices of irregularity in the rational
case.  We briefly investigated this issue in~\cite{Holden98}, where
there are tables of the analogue of $u_{r}(x)$ (using the index of
$\chi$-irregularity) for $x=1000$, $r$ from 0 to 4, and $D=5, 8, 12,$
and $13$.  The chi-squared test results are not included, but using
the methodology discussed earlier they are 3.32, 1.74, 1.15, and 2.54.
The corresponding significance levels are .345, .628, .765, and .469,
respectively.  We could total the values of (the analogue of)
$u_{r}(x)$ for the four values of $D$ and compare them to the
predicted values; we might expect that this would give us a better
significance level because of the larger ``sample size''.  However, in
this case the chi-squared result is 3.53 and the significance level is
.316, which is worse than any of the results for the values of $D$
taken separately!  This may be due to some small second-order
bias which is common to each sample and thus is reinforced when they
are pooled together.

\section{The Hypotheses:  New Results}
\label{newdata}
In the course of testing the algorithms in~\cite{mathcomp}, we
collected more data in addition to that above. Table~\ref{D5} shows
the number of primes less than 5000 which have $\chi$-index of
irregularity $r$ for various values of $r$ and $D=5$.  We compared the
observed and predicted distributions, using the methodology above,  for
primes below $x$ where $x$ was 1000, 2000, 3000, 4000, and 5000, and
found chi-squared values of 3.32, 5.03, 2.51, 1.73, and 2.10 and
significance levels of .344, .170, .473, .630, and .552, respectively.

\begin{table}[h!]
\caption{Results for $D=5$ and $p<5000$} \label{D5}
\begin{center}
\begin{tabular}{|c|c|c|c|}
\hline
$r$ & number & predicted number & predicted fraction  \\
\hline
0 & 422 & 405.16 &.606531  \\
1 & 186 & 202.58 &.303265  \\
2 & 51 & 50.65 &.075816  \\
3 & 7 & 8.44 &.012636  \\
4 & 2 & 1.06 &.001580  \\
\hline
\end{tabular}
\end{center}
\end{table}

Other data was obtained using the philosophy, described in
Section~\ref{algorithms}, of
computing the values of $L(1-2m, \chi)$ for large numbers of $D$ and
relatively small values of $m$.  As in the discussion of $D=5, 8, 12,$ and
$13$ above, we present the total across the different
discriminants.
Table~\ref{aveD5000} presents the data for all $D<5000$ and $p<100$.  The
chi-squared value for the totals is 81.1 and the significance level
is .000.

\begin{table}[h!]
\caption{Results for $D<5000$ and $p<100$} \label{aveD5000}
\begin{center}
\begin{tabular}{|c|c|c|c|}
\hline
$r$ & total & predicted total  & predicted  \\
& number & number               &        fraction\\
\hline
0 & 21864 & 22068.01 &.606531  \\
1 & 11596 & 11034.01 &.303265  \\
2 & 2529 & 2758.50  &.075816  \\
3 &   347 & 459.75 &.012636  \\
4 & 41 & 57.47  &.001580  \\
5 & 7 & 5.75 &   .000158 \\
\hline
\end{tabular}
\end{center}
\end{table}

However, if we view the data broken down by prime, as in
Table~\ref{breakdownD5000}, we see that a large part of the
contribution to the chi-squared value is from small primes.  The
values of $p$ shown in the table were selected with an eye towards
showing both a trend toward smaller chi-squared values as $p$ increases
and also some of the exceptions.  We hope to make the nature of the
small-prime contribution clearer in the future.

\begin{table}[h!]
\caption{Results for $D<5000$; selected values of $p$} \label{breakdownD5000}
\begin{center}
\begin{tabular}{|c|c|c|c|c|c|}
\hline
$r$ & 0 & 1 & 2 & $\geq 3$ & sig. level\\
\hline
pred.& 919.50 & 459.75 & 114.94 & 21.81& \\
\hline
$p=3$ & 876 & 640 & 0 & 0 & .000 \\
$p=5$ & 956 & 500 & 60 & 0 & .000 \\
$p=7$ & 895 & 530 & 89 & 2 & .000 \\
$p=11$ & 876 & 497 & 131 & 12 &.008\\
$p=13$ & 947& 467& 91 & 11&.010\\
$p=17$ &950& 452& 95& 19&.175\\
$p=23$&933& 462& 106&15& .387\\
$p=37$&913 & 468 & 108 & 27&.605\\
$p=47$&911 & 476 & 109 & 20&.775\\
$p=67$&915 & 466 & 114 & 21&.986\\
$p=79$ &859 & 487 & 144 & 26&.003\\
$p=97$ &909 & 468 & 122 & 17&.623\\
\hline
\end{tabular}
\end{center}
\end{table}

\section{Practical Notes}
\label{probabilities}
\label{practical}
The hypotheses that $p$ divides $\zeta(1-2m)$ with probability $1/p$
and that $p$ divides $L(1-2m, \chi)$ with the same probability clearly
imply that $p$ divides $\zeta_{D}(1-2m)= \zeta(1-2m) L(1-2m, \chi)$
with probability $(2/p)-(1/p^{2})$, or approximately $2/p$ for very
large $p$, as we claimed in Section~\ref{algorithms}.  Also, for the
algorithms in that section one
doesn't really have to worry about the possibility that $\delta$ is
not $p-1$, since this would require $D=p$.  However, we showed that
$D_{2}-D_{1}$ can be on the order of magnitude of $\log p$, so the
case of $D=p$ can only arise as the result of what can only be called
bad planning.

As mentioned in the introduction, one use for the algorithms of this
paper may be to find $D$ and $p$ such that the class group of
$\Q(\sqrt{D}, \zeta_{p})$ can be used for cryptographic protocols.
In~\cite{BP97}, Buchmann and Paulus introduced a one way function
based on class groups of number fields and noted that such a function
could be used to implement Diffie-Hellman key exchanges and ElGamal
signature schemes, to take two examples.  These ideas are expanded on
in~\cite{Meyer-roots}, which introduces a signature scheme called RDSA
which is based on taking $p$-th roots in the class group of a number
field or other abelian group.  Here $p$ is a random prime number which
(one assumes) does not divide the order of the group.  One advantage
of this signature scheme is that it is unnecessary and in fact
undesirable to know precisely the order of the group, a situation
which frequently occurs with class groups and in fact is generally
true for the class groups found with the algorithms above.

Given that the order of the class group is unknown, the question of
which class groups are suitable for these protocols is addressed
in~\cite{Meyer-fields}, which gives two necessary conditions on the
class number:
\begin{itemize}
    \item The class number must be sufficiently large.  This should
    make it difficult to determine the class number or discrete
    logarithms using exhaustive search, Pollard Rho,
    Baby-Step-Giant-Step, Hafner-McCurley, or index calculus.

    \item
    The class number must have at least one sufficiently large prime divisor.
    This should make it difficult to find discrete logarithms using a
    Pohlig-Hellman attack.
\end{itemize}
As we have seen, the algorithms of this paper allow us to find a
class number with a prime divisor as large as desired, and thus with
the class number itself as large as desired.

The drawback is that the amount of time and space needed to carry out
the cryptographic protocols in these groups can also be very large.
The papers~\cite{BP97} and~\cite{Meyer-fields} explain how to
represent the objects necessary to compute with.  Since elements in
the class group are equivalence classes of ideals in a ring of
integers, we need to store a $\Z$-basis for the ring of integers of
$\Q(\sqrt{D}, \zeta_{p})$. As noted in~\cite{BP97} and explained in
more detail in~\cite{BvS}, this requires
$(\log \left|\Delta\right|)^{O(1)}$ bits of storage, where $\Delta$ is
the discriminant of $\Q(\sqrt{D}, \zeta_{p})$.  Unfortunately, it is
not hard to show that $\Delta=D^{p-1}p^{2p-4}$ if $p$ does not divide
$D$ or $D^{p-1}p^{p-3}$ if $p$ does divide $D$ (but $p\neq D$).
(See~\cite{Washington}, for instance.) Thus
the $\Z$-basis requires $(p \log D)^{O(1)}$ bits of storage.
Furthermore, as explained in~\cite{Meyer-fields}, an ideal class
should be represented by a member of the class which is LLL-reduced;
that is, by one which corresponds to an LLL-reduced lattice under
Minkowski's embedding.  Such an representation requires
$(n+\log \left|\Delta\right|)^{O(1)}$ bits of storage, where $n$ is
the degree of the field and $\Delta$ is as before.  (See~\cite{BvS}
and~\cite{BK} for details.)  In our case $n=p-1$ (assuming
$p\neq D$) so this again requires $(p \log D)^{O(1)}$ bits of storage.
Of course, this means that the time it takes to carry out the basic
algorithms for the class group is also generally going to be
exponential in the size of $p$.  Whether this situation is bad enough
to preclude the use of our fields is not yet clear.

\section{Conclusion and Future Work}
\label{conclusion}
Much of the future work described in \cite{mathcomp} still remains to
be done; in particular many improvements could be made in the
implementations of the algorithms and perhaps in the algorithms
themselves.  However, the results already seem encouraging.  With a
faster implementation, the use of the search algorithm of
Section~\ref{algorithms} to find class groups large enough for secure
cryptography seems quite feasible, although this should be tested in
practice.  More importantly, an implementation of one or more
cryptographic protocols needs to be done using the class groups we
have described in order to determine whether secure cryptography can
be done sufficiently quickly in these groups.

The data collected in Section~\ref{newdata} is also encouraging, but
clearly more is necessary.  The author hopes to implement and run his
algorithms on a true supercomputer in the near future.  The data
produced by this will undoubtably give a clearer picture of the
phenomena so far observed, perhaps leading to refinements of our
hypotheses.

\section*{Acknowledgements}
The author would like to thank Don Burdick of the Institute of
Statistics and Decision Sciences at Duke University for his help in
making sense of the data presented in Sections~\ref{olddata}
and~\ref{newdata}.  He would also like to thank Carl Pomerance for
suggesting that this data was worthy of presenting to a statistician
in the first place.

\newcommand{\SortNoop}[1]{}
\providecommand{\bysame}{\leavevmode\hbox to3em{\hrulefill}\thinspace}

\end{document}